\documentclass[11pt,a4paper]{article}

\usepackage{amsfonts}
\usepackage{graphics}

\begin{document}

\title{EXACT MINIMAX RISK FOR DENSITY ESTIMATORS IN NON-INTEGER SOBOLEV CLASSES}
\author{\textbf {Clementine Dalelane}\\ \small{Laboratoire de Probabilit\'es
et Mod\`eles Al\'eatoires}\\
\small{Universit\'e Pierre et Marie Curie, Paris VI}\\
\small{dalelane@ccr.jussieu.fr}}
\date{\today}
\maketitle

\begin{abstract}\noindent The $L_2$-minimax risk in Sobolev classes
of densities with non-integer smoothness index is shown to have an
analog form to that in integer Sobolev classes. To this end, the
notion of Sobolev classes is generalized to fractional derivatives
of order $\beta\in\mathbb R^+$. A minimax kernel 
density estimator for such a classes is found. 
Although   
there exists no corresponding proof in the literature so far, the result
of this article was used implicitly in numerous papers. 
A certain necessity that this gap had to be filled, can thus not be
denied.\end{abstract}



\section{Introduction}

When trying to describe the goodness of an estimator, minimax
performance is one optimality criterium possible to be consulted.
The minimax risk of density estimators can be regarded in various
settings, e.g. we differentiate between the local risk in a single
point and the integrated risk over the whole curve. Several loss
functions have been under consideration, such as absolute,
quadratic and supremum norm, Hellinger and Kullback-Leibler
distance. But \underline{exact} asymptotics is up to now limited
to a few special cases: to supremum risk in H\"older classes, and
to mean integrated square error (MISE) in analytical and in
Sobolev classes.

The latter has been examined for quite a while since in 1983,
Efroimovich and Pinsker completed the asymptotic minimax rate of
the lower bound of MISE (Samarow \cite{sama}) by the still lacking
asymptotically exact constant, using tools that are common in
information theory. The results were enhanced and new methods of
proof found by Golubev \cite{gol1} and \cite{gol2}, Golubev, Levit
\cite{gol3} and Schipper \cite{sch}. Sobolev classes are classes
of $L_2$-integrable functions, in the present problem densities,
for which smoothness is measured through the $L_2$-norm of their
$\smash{\beta^{\mbox{th}}}$ derivative, $\beta\in\mathbb N$.
\begin{eqnarray}
\label{1} \mathcal S_{\beta}(L)&=&\left\{f\in L_2\Bigl|\int
\left(f^{(\beta)}(x)\right)^2dx
    \le L\right\},\quad L<\infty
\end{eqnarray}
Nowadays it is well known that
\begin{eqnarray}
    \label{2}
\inf_{\widetilde f_n}\sup_{f\in\mathcal
    S_{\beta}(L)}n^{\frac{2\beta}{2\beta+1}}
    \ E_f\|\widetilde f_n-f\|_2^2&=&\gamma(\beta,L)\Bigl(1+o(1)\Bigr)\\
    \nonumber
\mbox{where}\qquad\qquad\qquad\gamma(\beta,L)&=&(2\beta+1)\left(
    \frac{\pi(2\beta+1)(\beta+1)}{\beta}\right)^{-\frac{
    2\beta}{2\beta+1}}L^{\frac{1}{2\beta+1}}
\end{eqnarray}
is Pinsker's constant. Estimators attaining minimax rates of
convergernce have been studied in abundance, e.g. kernel
estimators, but also wavelet estimators and a wide range of
others. More care has to be taken when envisaging asymptotically
exact minimax estimators.

However, the characterization of the smoothness of a given
density function is incomplete when just assigning it to some
$\mathcal S_{\beta}(L)$, $\beta\in\mathbb N$. Recalling the
Sobolev criterion,
\begin{eqnarray*}
L&\ge&\int \left(f^{(\beta)}(x)\right)^2dx
    \ =\ \frac{1}{2\pi}\int |\omega^{\beta}\widehat f(\omega)|^2d\omega,
\end{eqnarray*}
we immediately observe that $\mathcal S_{\beta}(L)$ contains
densities which do not lie in $\mathcal S_{\beta+1}$, although for
suitably chosen $L'<\infty$ and $\varepsilon<1$, they certainly do
satisfy
$\smash{\frac{1}{2\pi}\int|\omega^{\beta+\varepsilon}\widehat
f(\omega)|^2d\omega\le L'}$.

The present article is interested in the question of whether the
minimax risk can also be calculated for such generalized Sobolev
classes. Corresponding claims are implicit in a number of recent
papers, yet their proofs cover but the entire case. For our
purpose we will employ the concept of the so-called fractional
derivative after Riemann and Liouville, thoroughly discussed in
Samko \cite{sam}:
\begin{eqnarray*}
f^{(\beta)}(x)&=&\frac{d^{\beta}}{dx^{\beta}}f(x)\ =\
    \frac{1}{\Gamma
    (\lceil\beta\rceil-\beta)}\ \frac{d^{\lceil\beta\rceil}}
    {dx^{\lceil\beta\rceil}}\int
    t^{\beta-\lceil\beta\rceil}f(x+t)dt
\end{eqnarray*}
with $\lceil x\rceil$ the smallest integer greater than the
positive real number $x$. For $\beta\in \mathbb N$, $f^{(\beta)}$
is the $\smash{\beta^{\mbox{th}}}$ derivative of $f$, for
$\beta\in\mathbb R^+\backslash\mathbb N$ it is the
$\smash{\beta^{\mbox{th}}}$ fractional derivative of $f$ (Samko
\cite{sam}, p. 137). In case $f^{(\beta)}$ is continuous and
$L_1$-integrable, then $\smash{\widehat
{f^{(\beta)}}}(\omega)=(-i\omega)^{\beta} \widehat f (\omega)$.
The other way around, if $(-i\cdot\mbox{id})^{\beta}\widehat f\ $
is $L_1$- or $L_2$-integrable, the inverse transform from the
Fourier into the time domain exists and for our purpose we define:
\begin{eqnarray}
    \label{3}
f^{(\beta)}(x)&:=&\frac{1}{2\pi}\int (-i\omega)^{\beta}\widehat f
    (\omega)e^{-ix\omega}dx
\end{eqnarray}
Existence and uniqueness of the $\smash{\beta^{\mbox{th}}}$
fractional derivative of $f$ follow thus from $\frac{1}{2\pi}\int
\smash{|\omega^{\beta}\widehat f(\omega)|^2}d\omega$ $\le L$, and
Parseval's equality gives $\smash{\int
(f^{(\beta)}(x))^2dx=\frac{1}{2\pi}\int |\omega^{\beta}\widehat
f(\omega)|^2d\omega}$.

Adopting the idea of Schipper \cite{sch}, we find upper and lower
bounds for the asymptotic minimax risk in $\mathcal
S_{\beta}(L),\ \beta
>1/2$, which are then shown to converge towards each other.
Thereby it will be verified that the minimax risk is determined by
$\smash{n^{2\beta/(2\beta+1)}\gamma(\beta,L)}$, where
$\smash{\gamma(\beta,L)}$ is an analogue of Pinsker's constant. A
minimax kernel function for kernel density estimation is obtained
as a byproduct from the calculation. On benefit of our a
statement, it is for instance possible to show the asymptotically
exact minimax-adaptivity of non-parametric estimation procedures
such as the recently proposed Stein's blockwise estimator for
densities (Rigollet \cite{rig}) and the cross-validation kernel
choice for density estimation (Dalelane \cite{dal}).

The calculation of the upper bound in Schipper \cite{sch} actually
holds for both entire and non-entire smoothness indeces, so Schipper's
Theorem 3 ({\sl Math. Meth. of Statistics} (1996), Vol.
5 No. 3, page 258-260) applies directly. To show the lower bound in
Section 2,
we replace original problem of estimating a curve by the
problem of estimating a finite-dimensional parameter $\theta$ (of
increasing dimension). A lower bound for the risk of such an
estimator may be found by means of the van Trees inequality. The
Bayesian risk over a least favorable parametric family of
densities $\mathcal F_{\Theta}$, and a least favorable prior
distribution $\Lambda$ on the space of finite-dimensional
parameters $\Theta$, such that $f_{\theta}\in \mathcal
S_{\beta}(L)$ with a high probability, provides us with a lower
bound for the minimax risk on $\mathcal S_{\beta}(L)$. It is
exactly this gap in the literature: $f_{\theta}$ asymptotically in
$\mathcal S_{\beta}(L)$ for $\beta\not\in\mathbb N$, which we have
been able to close in the present paper. Although the
demonstrations follow in general the same lines as Schipper
\cite{sch}, the least favorable family of densities had to be
constructed in a different way. The proof of the essential
property (Theorem 2) applies Riemann-Liouville calculus along with
approximations in the Fourier domain and is not similar to
Schipper \cite{sch}.

The result for the lower bound can be considered as a special case of
the theorem in Golubev \cite{gol2}, who yields lower bounds for
the quadratic risk of non-parametric estimation problems in a
variety of elliptic density classes via Local Asymptotic
Normality. Unfortunately the proof in Golubev \cite{gol2} is
heavily abbreviated (the proof of a claim corresponding to our
Theorem 2 is actually omitted) and not easy to retrace. We hope
that by our detailed proof, we are able to somehow enlighten the
complicated matters.

\section{Minimax bounds}

Let $X_1,\ldots,X_n$ be i.i.d. random variables with common
density function $f$ and let $\widetilde f_n$ be an arbitrary
estimator for $f$ depending but on the sample.
\\\\
{\sl\textbf{Theorem 1} (see Schipper \cite{sch} Theorem 3)
\quad Let $\mathcal S_{\beta}(L)$ be the
Sobolev class of those $L_2$-integrable densities, which satisfy
$\frac{1}{2\pi}\int |\omega^{\beta}\widehat f(\omega)|^2d\omega\le
L$ for some constants $\beta>0$ and $L<\infty$. Then it holds,
that}
\begin{eqnarray*}
\inf_{\widetilde f_n}\sup_{{\mathcal
S}_{\beta}(L)}n^{\frac{2\beta}{2\beta+1}}\
    E_f\|\widetilde f_n-f\|^2_2&\le&\gamma(\beta,L)
\end{eqnarray*}
{\sl The bound is maintained by a kernel estimator with the
minimax kernel $K_{\beta}$, that is the inverse Fourier transform
of $\widehat
K_{\beta}(\omega)=\left(1-c_{\min}(L)\cdot|\omega|^{\beta}\right)_+$
with $c_{\min}=\left(\frac{nL\pi(2\beta+1)(\beta+1)}{\beta}\right)
^{-\beta/(2\beta+1)}$.}
\\\\
Generally speaking, the derivation of the lower bound proceeds
similarly to Schipper \cite{sch} (Subsection 4.1, page 262-268).
However it is not the same, and so we give a little more
detail. The following steps lead to the
desired result, which can partly be effected analogously to
Schipper \cite{sch}, partly new proofs had to found:

1) Construction of a least favorable parametric family of
densities $\mathcal F_{\Theta}$, proof that the elements of
$\mathcal F_{\Theta}$ are contained in an
$\varepsilon$-neighborhood of the considered Sobolev class
$\mathcal S_{\beta}(L)$. Both our center function of $\mathcal
F_{\Theta}$ and our perturbation functions had to be constructed
in a distinct way to Schipper \cite{sch}, whereas the parameter
set $\Theta$ is the same. The proof of our Theorem 2 is different
to that of Schipper's corresponding lemmata (Lemma 1 through 4).

2) Definition of a least favorable prior distribution $\Lambda$ on
the parameter set $\Theta$, proof that under $\Lambda$ the
elements of $\mathcal F_{\Theta}$ are contained in $\mathcal
S_{\beta}(L)$ itself with high probability. This time, Schipper's
distribution $\Lambda$ and his proof (Lemma 5, p. 266-267) are
possible to transfer to our context.

3) Main approximation of the lower bound via the Bayes risk over
$\mathcal F_{\Theta}$ with respect to $\Lambda$ by means of the
van Trees inequality. Again the proof of Schipper's Proposition 2
(p. 267-268) resembles our demonstration.
\\\\
The problem of searching a lower bound for the minimax risk over
the Sobolev class $\mathcal S_{\beta}(L)$, can be reverted to a
parametric subset of $\mathcal S_{\beta}$. Whether the minimax
risk over the subclass coincides with the minimax risk over
$\mathcal S_{\beta}(L)$, obviously depends on the difficulty of
the estimation problem within the subclass. We achieve our aim
using the adjacent construction:

Let us assume $\beta>1/2$ and let $f_0$ be the following density
from $\mathcal S:=\bigcap\limits_{\beta\in\mathbb R^+}\mathcal
S_{\beta}$
\begin{eqnarray}
    \label{5}
f_0(x)&:=&\left\{\begin{array}{l}\frac{1}{c_a}\exp\Bigl\{-\frac{a}
    {(x+1/2)(1/2-x)}\Bigr\},\ -1/2\le x\le 1/2\\
        0,\ \mbox{otherwise}
    \end{array}\right.
\end{eqnarray}
with $c_a$ defined such that $f_0$ is a density, and for technical
reason, the constant $a$ satisfying $\int |\widehat
f_0(\omega)\frac{\sin\omega/2}{\omega/2}| d\omega=2\pi$. Since
$\int |\omega^{\beta}\widehat f_0(\omega)|^2d\omega<\infty$ for
all $\beta<\infty$, also $\int |\omega^{\beta}\widehat
f_0(\omega)|d\omega$ exists for all $\beta<\infty$. Let $g_A$ be
the indicator function on $[-A+1/2,A-1/2]$ times the factor
$\frac{1}{2A-1}$, i.e.
\begin{eqnarray}
    \label{6}
g_A(x)&:=&\frac{1}{2A-1}\ I_{\left[-A+\frac{1}{2}\,,\,
    A-\frac{1}{2}\right]}(x).
\end{eqnarray}
Then $f_0\ast g_A$ is a symmetric density within $\mathcal S$,
that takes the constant value $\smash{\frac{1}{2A-1}}$ on
 $[-A+1,A-1]$, and decreases smoothly towards 0 on  $[A-1,A]$ and
$[-A,-A+1]$. In order to constitute a sufficiently difficult
estimation problem departing from this very smooth density, let us
add some perturbation functions to $f_0\ast g_A$:
\begin{eqnarray}
    \label{7}
\varphi_{k}(x)&:=&\left\{\begin{array}{ll}\frac{1}{\sqrt{A}}\cos\frac{k\pi x}{A}\ I_{[-A,A]}(x),& k>0\\
        \frac{1}{\sqrt{A}}\sin\frac{k\pi x}{A}\ I_{[-A,A]}(x),& k<0\end{array}\right.
\end{eqnarray}
These perturbations will be weighted by factors $\theta_k$, where
$\theta=(\ldots,\theta_{-2},\theta_{-1},\theta_1,\theta_2,\ldots)$
is (asymptotically) in the set:
\begin{eqnarray}
    \label{8}
\Theta_A(L)&:=&\left\{\theta\in \mathbb R^{\infty}\Bigl|\
\sum_{k\neq0}|\theta_k|\le
    A^{-2\beta+1}\mbox{ and }
    \sum_{k\neq0}\theta_k^2\left(\frac{k\pi}{A}\right)^{2\beta}
    \le 4A^2L\right\}
\end{eqnarray}
The set $\{f_{\theta}|\,\theta\in \Theta_A(L)\}$ will from now on
be the family of densities under consideration.
\begin{eqnarray}
    \label{9}
f_{\theta}(x)&:=&\frac{1}{b(\theta)}\ f_0\ast
g_A(x)\Bigl(1+\sum\limits_{k\neq0} \theta_k\varphi_k(x)\Bigr),
\end{eqnarray}
where $b(\theta)$ is the normalizing constant. We
\underline{cannot} prove that
 $\{f_{\theta}|\,\theta\in\Theta_A(L)\}\subseteq \mathcal S_{\beta}(L)$, but instead
that for all $\varepsilon>0$ there exists an
$A_{\varepsilon}<\infty$, so that for every $A\ge A_{\varepsilon}$
the following holds:
$\sup_{\theta\in\Theta_A(L)}\|f_{\theta}^{(\beta)}\|_2^2\le
L+\varepsilon$.
\\ \\
{\sl\textbf{Theorem 2}\quad Let $f_{\theta}$ and $\Theta_A(L)$ be
defined as above. Then, as $A\longrightarrow\infty$:}
\begin{eqnarray*}
\sup_{\Theta_A(L)}\|f_{\theta}^{(\beta)}\|_2^2&=&L+o(1)
\end{eqnarray*}
This theorem is the main assertion of our paper. Filling the gap
in the hitherto existing literature, it enables us to go on
proving the minimax bound for non-integer Sobolev classes. Its
cumbersome and unpleasantly lengthy proof is to be found in
Section 5.

The next step leading to the lower bound requires the definition
of a prior distribution $\Lambda$, which is done accordingly to
\cite{sch}, so as to yield a parameter $\theta$ of finite
dimension: Let $\varepsilon>0$, $W>0$ and $\sigma_k^2>0$
\begin{eqnarray}
    \label{10}
\lambda(\theta)&=&\prod\limits_{0<|k|<W}\lambda_k(\theta_k)\
\prod\limits_{|k|\ge W}\delta_0(\theta_k),
\end{eqnarray}
where $\delta_0(.)$ is the Dirac function on 0, and for $|k|<W$:
$\lambda_k(\theta_k)$ are absolutely continuous densities with
$E\theta_k^2=\sigma_k^2$, $\theta_k^2\le G^2\sigma_k^2$
($\Lambda$-f.s.) for some $G<\infty$, and the Fisher information
$\smash{I_k:=\int\frac{{\lambda'_k}^2(\theta_k)}{\lambda_k(\theta_k)}
\,d\theta_k}$ $\smash{\le(1+\varepsilon)\sigma_k^{-2}}$ (with
respect to the translation group $\{\lambda_k(. -u)|u\in\mathbb
R\}$). (These conditions are satisfied, for example, by
independent bounded, zero mean random variables $\sigma_k\xi_k$,
$|k|<W$, with $|\xi_k|<G$, $E\xi_k^2=1$ and the Fisher-information
of the density of $\xi_k$ smaller than $1+\varepsilon$.) Let us
set
\begin{eqnarray}\nonumber
W&=&\frac{A}{\pi}\Bigl(\frac{L(1-\varepsilon)n(2\beta+1)(\beta+1)\pi}{\beta}\Bigr)^
    {\frac{1}{2\beta+1}}\\
    \label{11}
&&\\\nonumber
\sigma_k^2&=&\frac{4A}{n}\Bigl(\Bigl|\frac{W}{k}\Bigr|^{\beta}-1\Bigr)_+\\\nonumber
\end{eqnarray}
As $W$ grows with $n\longrightarrow\infty$, the dimension of the
parameter $\theta$ will tend to infinity, allowing for more and
more perturbation functions $\varphi_k$ in the definition of
$f_{\theta}$. At the end of Section 3 it will be shown that
$\sigma_k^2$ and $W$ of this form approximately maximize the lower
bound of the minimax risk for the prior distribution $\Lambda$.

Since $\Lambda$ is not supported on $\Theta_A(L)$, we will have to
show that at least the probability of $\theta\in\Theta_A(L)$ grows
with $n\longrightarrow\infty$:

First consider that $\lambda$ has a bounded support, $|\theta_k|
\le G\sigma_k$ for $|k|<W$, and else $\theta_k=0$. With the above
construction of $\sigma_k^2$ and $W$, letting $A\sim\ln n$,
condition $\smash{\sum|\theta_k|\le A^{-2\beta+1}}$ is fulfilled
for $n$ sufficiently large. Lemma 1 takes care of
$\smash{\sum\theta_k^2
\left(\frac{k\pi}{A}\right)^{2\beta}\le4A^2L}$.
\\ \\
{\sl\textbf{Lemma 1}\quad For the prior distribution $\Lambda$
defined above, with $W$ and $\sigma_k^2$ as in (10), it holds that
for $n\longrightarrow\infty$:}
\begin{eqnarray*}
P_{\lambda}\Bigl(\theta\not\in \Theta_A(L)\Bigr)&=&o(n^{-1})
\end{eqnarray*}
This lemma corresponds to Lemma 5 in Schipper \cite{sch}, p. 266.
Its proof is exactly the same (p. 266-267) and we abstain from
quoting it here (also see Dalelane \cite{dal} for more details).
\\\\
{\sl\textbf{Theorem 3}\quad For $L<\infty$, $\beta>1/2$ and
$\gamma(\beta,L)$ equal to Pinsker's constant we have:}
\begin{eqnarray*}
\liminf_{n\rightarrow \infty}\ \inf_{\widetilde
f_n}\sup_{f\in\mathcal S_{\beta}(L)}
    n^{\frac{2\beta}{2\beta+1}}\ E_f\|\widetilde f_n-f\|_2^2
    &\ge&\gamma(\beta,L)
\end{eqnarray*}
\textbf{Proof}\quad Let us at first reduce the supremum of the
risk by restricting the set of density functions. According to
Theorem 2 we know that for $A\sim\ln n$,
$\lim_{A\rightarrow\infty}\{
f_{\theta}|\theta\in\Theta_A(L)\}\subseteq\mathcal S_{\beta}(L)$.
\begin{eqnarray}
    \label{12}
    \liminf_{n\rightarrow\infty}\inf_{\widetilde
    f_n}\sup_{S_{\beta}(L)}E_f\|\widetilde f_n-f\|_2^2
    &\ge&\liminf_{n\rightarrow\infty}\inf_{\widetilde f_n}\sup_{\Theta_A(L)}
    E_{f_{\theta}}\|\widetilde f_n-f_{\theta}\|_2^2
\end{eqnarray}
For any fixed $A$, we find a lower bound for the supremum over
$\Theta_A(L)$ through the Bayesian risk with respect to $\Lambda$.
\begin{eqnarray}
    \label{13}
\inf_{\widetilde f_n}\sup_{\Theta_A(L)}
    E_{f_{\theta}}\|\widetilde f_n-f_{\theta}\|_2^2
&\ge&\inf_{\widetilde
    f_n}\int_{\Theta_A(L)}E_{f_{\theta}}\|\widetilde
    f_n-f_{\theta}\|_2^2\ d\Lambda(\theta)
\end{eqnarray}
In inequality (\ref{21}) of the proof of Theorem 2 it will be
shown that $1-A^{-3/2}\le b(\theta)\le 1+A^{-3/2}$. Furthermore, because of
orthonormality
$\|\sum\theta_k\varphi_k\|_2^2=\sum\theta_k^2\le\sum|\theta_k|\le
A^{-2\beta+1}$, (\ref{34}). So we can derive, for all
$\theta\in\Theta_A(L)$:
\begin{eqnarray*}\nonumber
\|f_{\theta}\|_2 &=&\frac{1}{b(\theta)}\Bigl\|f_0\ast
    g_A\Bigl(1-\sum_{k\neq
    0}\theta_k\varphi_k\Bigr)\Bigr\|_2\\\nonumber
&\le&\frac{1}{b(\theta)}\ \max f_0\ast g_A\ \Bigl\|I_{[-A,A]}+
    \sum_{k\neq
    0}\theta_k\varphi_k\Bigr\|_2\\\nonumber
&\le&\frac{\mbox{const.}}{\sqrt{A}}\ =:\ \frac{1}{\sqrt{A_0}}
\end{eqnarray*}
Because the set of all densities with $\|f\|_2\le1/A_0$ is convex,
we may in (\ref{14}) also restrict the set estimators to
$\|\widetilde f_n\|_2^2\le1/A_0$ without increasing the supremum.
\begin{eqnarray}
    \label{14}
\hspace{-.3cm}\inf_{\widetilde
f_n}\int_{\Theta_A(L)}E_{f_{\theta}}\|\widetilde
    f_n-f_{\theta}\|_2^2\ d\Lambda(\theta) &=&\inf_{\|\widetilde
    f_n\|_2^2\le
    A_0^{-1}}\int_{\Theta_A(L)}E_{f_{\theta}}\|\widetilde
    f_n-f_{\theta}\|_2^2\
    d\Lambda(\theta)\qquad\qquad\qquad\qquad\qquad\quad\\
    \label{15}
&\ge&\inf_{\|\widetilde f_n\|_2^2\le A_0^{-1}}\int
    E_{f_{\theta}}\|\widetilde f_n-f_{\theta}\|_2^2\
    d\Lambda(\theta) - \frac{4}{A_0}\,P_{\lambda}\Bigl(\theta\not\in\Theta_A(L)     \Bigr)\\\nonumber
&=&\inf_{\|\widetilde f_n\|_2^2\le A_0^{-1}}\int
    E_{f_{\theta}}\|\widetilde f_n-f_{\theta}\|_2^2\
    d\Lambda(\theta)\ +\ o(n^{-1})\\
    \label{16}
&\ge&\inf_{\widetilde f_n}E_{\lambda} E_{f_{\theta}}\|\widetilde
    f_n-f_{\theta}\|_2^2
    \ +\ o(n^{-1})
\end{eqnarray}
Due to $\|f_{\theta}\|_2^2\le A_0^{-1}$ and $\|\widetilde
f_n\|_2^2\le A_0^{-1}$ it holds in (\ref{15}) that $\|\widetilde
f_n-f_{\theta}\|_2^2\le 4A^{-1}_0$. In (\ref{16}) we return to the
complete set of estimators.

Since $f_{\theta}$ has bounded support, i.e. $[-A,A]$, it is
equivalent, as regards the quadratic risk, either to estimate the
function $f_{\theta}$ in the time domain or its Fourier
coefficients. ($\widehat f_{\theta}(0)=1$ is known)
\begin{eqnarray}\nonumber
&&E_{\lambda} E_{f_{\theta}}\|\widetilde
    f_n-f_{\theta}\|_2^2\\\nonumber &=&E_{\lambda} E_{f_{\theta}}\
    \frac{1}{2A}\sum_{\kappa\neq0}
    \Bigl|\widetilde{\widehat f_n}\left(\frac{\kappa\pi}{A}\right)-\widehat f_{\theta}\left(\frac{\kappa\pi}
    {A}\right)\Bigr|^2\\\nonumber
&=&E_{\lambda} E_{f_{\theta}}\ \frac{1}{2A}\sum_{\kappa\neq0}\
    \mbox{Re}^2\left(\widetilde{\widehat f_n}
    \left(\frac{\kappa\pi}{A}\right)-\widehat f_{\theta}\left(\frac{\kappa\pi}
    {A}\right)\!\right)+\ \mbox{Im}^2\left(\widetilde{\widehat f_n}
    \left(\frac{\kappa\pi}{A}\right)-\widehat f_{\theta}\left(\frac{\kappa\pi}
    {A}\right)\!\right)\\
    \label{17}
&=&E_{\lambda} E_{f_{\theta}}\
    \frac{1}{2A}\sum_{\kappa\neq0}\left(\mbox{Re}\widetilde{\widehat
    f_n}
    \left(\frac{\kappa\pi}{A}\right)-\mbox{Re}\widehat f_{\theta}\left(\frac{\kappa\pi}
    {A}\right)\!\right)^2+\ \left(\mbox{Im}\widetilde{\widehat f_n}
    \left(\frac{\kappa\pi}{A}\right)-\mbox{Im}\widehat f_{\theta}\left(\frac{\kappa\pi}
    {A}\right)\!\right)^2\qquad\qquad
\end{eqnarray}
The van Trees inequality (Gill, Levit \cite{gil}) may now be
applied on every single summand. For technical reason, the real
parts are derived with respect to $\theta_{|\kappa|}$, while the
imaginary ones are derived with respect to $\theta_{-|\kappa|}$.
\begin{eqnarray}
    \nonumber
E_{\lambda} E_{f_{\theta}}
    \left[\mbox{Re}\widetilde{\widehat f_n}\left(\frac{\kappa\pi}{A}\right)-\mbox{Re}\widehat f_{\theta}
    \left(\frac{\kappa\pi}
    {A}\right)\right]^2&\ge&\frac{E_{\lambda}^2\left[\partial\,\mbox{Re}\widehat f_{\theta}
    \left(\frac{\kappa\pi}
    {A}\right)/\partial \theta_{|\kappa|}\right]}{nE_{\lambda}I_{f_{\theta}}(\theta_{|\kappa|})
    +I_{|\kappa|}}\\
    \label{18}
&&\\
    \nonumber
E_{\lambda} E_{f_{\theta}}
    \left[\mbox{Im}\widetilde{\widehat f_n}\left(\frac{\kappa\pi}{A}\right)-\mbox{Im}\widehat f_{\theta}
    \left(\frac{\kappa\pi}
    {A}\right)\right]^2&\ge&\frac{E_{\lambda}^2\left[\partial\,\mbox{Im}\widehat f_{\theta}
    \left(\frac{\kappa\pi}
    {A}\right)/\partial
    \theta_{-|\kappa|}\right]}{nE_{\lambda}I_{f_{\theta}}(\theta_{-|\kappa|})+I_{-|\kappa|}}
\end{eqnarray}
where we denote
$\smash{I_{f_{\theta}}(\theta_{\kappa})=\int\frac{(\partial
f_{\theta}(x)/\partial \theta_{\kappa})^2} {f_{\theta}(x)}\,dx}$.
$I_{\kappa}$ is the ``Fisher information'' of $\lambda_{\kappa}$
and by construction $\le (1+\varepsilon)\sigma_{\kappa}^{-2}$ for
$|{\kappa}|<W$ and $=\infty$ for $|{\kappa}| \ge W$, respectively.
Hence all summands with $|{\kappa}|\ge W$ vanish from the sum.
Approximations for
 $I_{f_{\theta}}(\theta_{\kappa})$,
$\partial\,\mbox{Re}\,\widehat
f_{\theta}\left(\frac{{\kappa}\pi}{A}\right) /\partial
\theta_{|\kappa|}$ und $\partial\,\mbox{Im}\,\widehat
f_{\theta}\left(\frac{{\kappa}\pi}{A}\right)/\partial
\theta_{-|\kappa|}$ are available from \\ \\
{\sl\textbf{Lemma 2}\quad For $A\longrightarrow\infty$:}
\begin{eqnarray}\nonumber
I_{f_{\theta}}(\theta_{\kappa}) \ =\
    \frac{1+o(1)}{2A}\qquad\qquad\quad
\frac{\partial\,\mbox{Re}\widehat
    f_{\theta}\left(\frac{{\kappa}\pi}{A}\right)}{\partial
    \theta_{|\kappa|}}
    &=&\frac{1+o(1)}{2\sqrt{A}}\qquad\qquad\quad
\frac{\partial\,\mbox{Im}\widehat
    f_{\theta}\left(\frac{{\kappa}\pi}{A}\right)}{\partial
    \theta_{-|\kappa|}} \ =\ \frac{1+o(1)}{2\sqrt{A}}
\end{eqnarray}
\\
with $o(1)$ independent of $\kappa$ and $\theta_{\kappa}$. The
proof is postponed to Section 5. From (\ref{16}) completed by
(\ref{17}), the van Trees approximation (\ref{18}) and Lemma 2 we
thus have:
\begin{eqnarray}\nonumber
\hspace{-.2cm}\inf_{\widetilde f_n} E_{\lambda}
E_{f_{\theta}}\|\widetilde
    f_n-f_{\theta}\|_2^2 &=&\inf_{\widetilde f_n} E_{\lambda}
    E_{f_{\theta}}\frac{1}{2A}\sum_{\kappa\neq0}
    \Bigl|\widetilde{\widehat f_n}\left(\frac{\kappa\pi}{A}\right)-\widehat f_{\theta}\left(\frac{\kappa\pi}
    {A}\right)\Bigr|^2\\\nonumber
&=&E_{\lambda} E_{f_{\theta}}\frac{1}{2A}\sum_{\kappa\neq0}\!
    \left(\!\mbox{Re}\widetilde{\widehat f_n}\!
    \left(\frac{\kappa\pi}{A}\right)\!-\mbox{Re}\widehat f_{\theta}\!\left(\frac{\kappa\pi}
    {A}\right)\!\right)^2\!+\!\left(\!\mbox{Im}\widetilde{\widehat f_n}\!
    \left(\frac{\kappa\pi}{A}\right)\!-\!\mbox{Im}\widehat f_{\theta}\!\left(\frac{\kappa\pi}
    {A}\right)\!\right)^2\\\nonumber
&\ge&\frac{1}{2A}\sum_{\kappa\neq0}\ \frac{E_{\lambda}^2
    \left[\partial\,\mbox{Re}\widehat f_{\theta}
    \left(\frac{\kappa\pi}
    {A}\right)/\partial \theta_{|\kappa|}\right]}{nE_{\lambda}I_{f_{\theta}}(\theta_{|\kappa|})+I_{|\kappa|}}
    +\frac{E_{\lambda}^2\left[\partial\,\mbox{Im}\widehat f_{\theta}\left(\frac{\kappa\pi}
    {A}\right)/\partial \theta_{-|\kappa|}\right]}{nE_{\lambda}I_{f_{\theta}}(\theta_{-|\kappa|})
    +I_{-|\kappa|}}\\\nonumber
&=&\frac{1}{2A}\sum_{0<|\kappa|<W}\
    \frac{\left(\frac{1+o(1)}{2\sqrt{A}}\right)^2}
    {n\ \frac{1+o(1)}{2A}+(1+\varepsilon)\sigma_{|\kappa|}^{-2}}
    +\frac{\left(\frac{1+o(1)}{2\sqrt{A}}\right)^2}
    {n\ \frac{1+o(1)}{2A}+(1+\varepsilon)\sigma_{-|\kappa|}^{-2}}\\
    \label{19}
&=&\frac{1+o(1)}{2A(1+\varepsilon)}\sum_{0<|\kappa|<W}\frac{1}{n+2A\sigma_{\kappa}^{-2}}\
\end{eqnarray}
All sums obtained from $W$ and $\sigma_{\kappa}^2$ through
(\ref{19}), i.e. from a prior distribution $\Lambda$ satisfying
Lemma 1, are thus lower bounds of the minimax risk.

What we are searching for is a bound as large as possible, we
hence maximize (\ref{19}) subject to the constraint
$\smash{\sum\sigma_{\kappa}^2\left(\frac{\kappa\pi}{A}\right)^{2\beta}\le(1-\varepsilon)4A^2L}$,
such that $P(\theta\not\in\Theta_A(L))=o(n^{-1})$ remains valid.
The solution to this problem is $W$ and $\sigma_{\kappa}^2$ from
(\ref{11}). The maximum in (\ref{19}) can be approximated as
follows:
\begin{eqnarray}\nonumber
&&\frac{1}{2A(1+\varepsilon)n}\sum_{0<|\kappa|<W}\left(\Bigl|\frac{W}{\kappa}\Bigr|^{\beta}-1\right)
    \Bigl|\frac{\kappa}{W}\Bigr|^{\beta}\\\nonumber
&=&\frac{1}{A(1+\varepsilon)n}\sum_{0<\kappa<W}\left(1-\Bigl(\frac{\kappa}{W}\Bigr)^{\beta}\right)
     \\\nonumber
&=&\frac{1}{A(1+\varepsilon)n}\ \frac{\beta}{\beta+1}\
    W\Bigl(1+o(1)\Bigr)\\\nonumber
&=&(2\beta+1)\left(\frac{(2\beta+1)(\beta+1)\pi}{\beta
    n}\right)^{-\frac{2\beta}{2\beta+1}}
    L^{\frac{1}{2\beta+1}}\ \frac{(1-\varepsilon)^{\frac{1}{2\beta+1}}}
    {1+\varepsilon}\Bigl(1+o(1)\Bigr)\\
&=&n^{-\frac{2\beta}{2\beta+1}}\
    \gamma(\beta,L)\ \Bigl(1+o(1)\Bigr)    \label{20}
\end{eqnarray}
Combining (\ref{12}) with (\ref{13}), (\ref{16}), (\ref{19}) and
(\ref{20}), we obtain the required result:
\begin{eqnarray*}
\qquad\qquad\qquad\liminf_{n\rightarrow \infty}\ \inf_{\widetilde
    f_n}\sup_{\mathcal S_{\beta}(L)}
    \ n^{\frac{2\beta}{2\beta+1}}\ E_f\|\widetilde f_n-f\|_2^2
    &\ge&\gamma(\beta,L)\qquad\qquad\qquad\qquad\quad\qquad\square
\end{eqnarray*}

\section{Remaining Proofs}

\textbf{Proof} of \textbf{Theorem 2}\quad For $f_{\theta}$ defined
in equation (\ref{9}), it holds that
\begin{eqnarray*}
\|f_{\theta}^{(\beta)}\|_2&=&\frac{1}{b(\theta)}\Bigl\|\left(f_0\ast
g_A\right)^{(\beta)}+
    \Bigl(f_0\ast g_A\sum\limits_{k\neq0} \theta_k\varphi_k\Bigr)^{(\beta)}\Bigr\|_2\\
&\le&\frac{1}{b(\theta)}\Bigl\|\left(f_0\ast
g_A\right)^{(\beta)}\Bigr\|_2+
    \frac{1}{b(\theta)}\Bigl\|\Bigl(f_0\ast g_A\sum\limits_{k\neq0} \theta_k\varphi_k\Bigr)
    ^{(\beta)}\Bigr\|_2\\
\end{eqnarray*}
$b(\theta)$, $\|(f_0\ast g_A)^{(\beta)}\|_2^2$ and $\|(f_0\ast
g_A\,\sum \theta_k\varphi_k )^{(\beta)}\|_2^2$ are then considered
one by one. Remember definition (\ref{7}):
$\varphi_k(x)=A^{-1/2}\cos (\pi k/A)I(|x|\le A)$ for $k>0$, and
the same with sine for $k<0$. Take first the normalizing constant
$b(\theta)$:
\begin{eqnarray}\nonumber
b(\theta)&=&\int_{-A}^A f_0\ast
g_A(x)\left(1+\sum_{k\neq0}\theta_k\varphi_k(x)\right)dx\\\nonumber
&=&1+\sum_{k\neq0}\theta_k\int_{-A}^A f_0\ast
g_A(x)\varphi_k(x)dx\\\nonumber &=&1+\sum_{k>0}\theta_k\int_{-A}^A
f_0\ast g_A(x)\frac{1}{\sqrt{A}}\,\cos\frac{k\pi x}{A}\ dx
    +\sum_{k<0}\theta_k\int_{-A}^A f_0\ast g_A(x)\frac{1}{\sqrt{A}}\,\sin\frac{k\pi x}{A}\ dx\\\nonumber
&=&1+\sum_{k>0}\theta_k\int_{-A}^A f_0\ast
g_A(x)\frac{1}{\sqrt{A}}\,\cos\frac{k\pi x}{A}\ dx+\,0\\\nonumber
&=&1+\frac{1}{\sqrt{A}}\sum_{k>0}\theta_k\left[\int_{-A}^A
\frac{1}{2A-1}\ \cos\frac{k\pi x}{A}\ dx-
    2\int_{A-1}^A \left(\frac{1}{2A-1}-f_0\ast g_A(x)\right)\cos\frac{k\pi x}{A}\ dx\right]\\\nonumber
&=&1+\frac{1}{\sqrt{A}}\sum_{k>0}\theta_k\left[\,0-
    2\int_{A-1}^A \left(\frac{1}{2A-1}-f_0\ast g_A(x)\right)\cos\frac{k\pi x}{A}\ dx\right]
\end{eqnarray}
For the second term on the right-hand side we have:
\begin{eqnarray}
    \nonumber
\frac{2}{\sqrt{A}}\Bigl|\sum_{k>0}\theta_k
    \int_{A-1}^A \left(\frac{1}{2A-1}-f_0\ast g_A(x)\right)\cos\frac{k\pi x}{A}\
    dx\Bigr|
&\le&\frac{2}{\sqrt{A}}\sum_{k>0}|\theta_k|\
    \int_{A-1}^A \frac{1}{2A-1}\Bigl|\cos\frac{k\pi x}{A}\Bigr|\ dx\\
        \nonumber
&\le&\frac{2}{\sqrt{A}(2A-1)}\sum_{k>0}|\theta_k|\\
    \nonumber
&\le&\frac{1}{\sqrt{A}(A-1/2)}\ A^{-2\beta+1},
\end{eqnarray}
so that for $\beta>1/2$ and $A$ sufficiently large, it follows
that
\begin{equation}
    \label{21}
1-A^{-3/2}\ \le\ b(\theta)\ \le\ 1+A^{-3/2}
\end{equation}
$(f_0\ast g_A)^{(\beta)}$ is integrable in $L_2$. So instead of
the $L_2$-norm of $f_0\ast g_A$ in the time domain, by Parseval's
equality we may as well study the $L_2$-norm of its Fourier
transform.
\begin{eqnarray}
\nonumber \left\|(f_0\ast
    g_A)^{(\beta)}\right\|_2^2&=&\int_{-A}^A\Bigl|(f_0\ast g_A)
    ^{(\beta)}(x)\Bigr|^2dx\\\nonumber
&=&\frac{1}{2\pi}\int\Bigl|\omega^{\beta}\widehat
    f_0(\omega)\widehat g_A(\omega)\Bigr|^2d\omega\\\nonumber
&=&\frac{1}{2\pi}\int\Bigl|\omega^{\beta}\widehat
    f_0(\omega) \frac{2\sin(A-1/2)\omega}{(2A-1)\omega}\Bigr|^2d\omega\\\nonumber
&=&\left(\frac{1}{A-1/2}\right)^2\frac{1}{2\pi}\int\Bigl|\omega^{\beta-1}\widehat
    f_0(\omega)\sin(A-1/2)\omega \Bigr|^2d\omega\\\nonumber
&\le&\left(\frac{1}{A-1/2}\right)^2\frac{1}{2\pi}\int\Bigl|\omega^{\beta-1}\widehat
    f_0(\omega)\Bigr|^2d\omega\\\label{36}
&\le&\left(\frac{1}{A-1/2}\right)^2\|f_0^{(\beta-1)}\|_2^2
\end{eqnarray}
For $\beta\ge1$, clearly $\|f_0^{(\beta-1)}\|_2^2<\infty$ because
$f_0$ lies in $\mathcal S$. For $1/2<\beta<1$ we can calculate
$\|f_0^{(\beta-1)}\|_2^2=\frac{1}{2\pi}\int|\omega^{\beta-1}\widehat
f_0(\omega)|^2d\omega\le\|f_0\|_2^2+\frac{1}{\pi}|2\beta-1|^{-1}$,
which is also less than infinity.

The consideration of the last and most important term $\|(f_0\ast
g_A\cdot\sum\theta_k\varphi_k)^{(\beta)}\|_2^2$ requires a little
knowledge about fractional derivatives. For two sufficiently
regular functions $f$ and $g$, the Leibnitz formula takes the
following form:
\begin{eqnarray}\nonumber
\left(f\cdot g\right)^{(\beta)}&=&\sum_{i=0}^{\infty}
{\beta\choose i}\ f^{(i)}\cdot
    g^{(\beta-i)}
\end{eqnarray}\nonumber
where $\beta\choose i$ an analogue to the binomial coefficient
with natural numbers:
\[{\beta\choose i}\ =\ \frac{\beta!}{i!(\beta-i)!}\ =\ \frac{\beta(\beta-1)
    (\beta-2)\cdots}{i!\ (\beta-i)(\beta-i-1)\cdots}
    \ =\ \frac{\beta\cdots (\beta-i+1)}{i!}\]
As usual, ${\beta\choose0}=1$. Now we apply this expansion to
$(f_0\ast g_A\cdot\sum\theta_k\varphi_k)^{(\beta)}$. Recall the
definition: $\lfloor x\rfloor$ is the integer part of a real
number $x$, and for $x$ positive (as in our case) $\lceil
x\rceil:=\lfloor x\rfloor+1$.
\begin{eqnarray}\nonumber
&&\Bigl\|\Bigl(f_0\ast g_A\cdot\sum_{k\neq0}
    \theta_k\varphi_k\Bigr)^{(\beta)}\Bigr\|_2\\\nonumber
&=&\Bigl\|\sum_{i=0}^{\infty} {\beta\choose i}
    \ (f_0\ast g_A)^{(i)}\Bigl(\sum_{k\neq0}
    \theta_k\varphi_k\Bigr)^{(\beta-i)}\Bigr\|_2\\\nonumber
&\le&\sum_{i=0}^{\lfloor\beta\rfloor} {\beta\choose i}
    \ \Bigl\|(f_0\ast g_A)^{(i)}\Bigl(\sum_{k\neq0}
    \theta_k\varphi_k\Bigr)^{(\beta-i)}\Bigr\|_2+
    \Bigl\|\sum_{i=\lceil\beta\rceil}^{\infty} {\beta\choose i}
    \ (f_0\ast g_A)^{(i)}\Bigl(\sum_{k\neq0}
    \theta_k\varphi_k\Bigr)^{(\beta-i)}\Bigr\|_2\\\nonumber
&\le&\sum_{i=0}^{\lfloor\beta\rfloor} {\beta\choose i}
    \ \max|(f_0\ast g_A)^{(i)}|\,\Bigl\|\Bigl(\sum_{k\neq0}
    \theta_k\varphi_k\Bigr)^{(\beta-i)}\Bigr\|_2+
    \Bigl\|\sum_{i=\lceil\beta\rceil}^{\infty} {\beta\choose i}
    \ (f_0\ast g_A)^{(i)}\Bigl(\sum_{k\neq0}
    \theta_k\varphi_k\Bigr)^{(\beta-i)}\Bigr\|_2\\\nonumber
&\le&\sum_{i=0}^{\lfloor\beta\rfloor} {\beta\choose i}
    \ \max|g_A|\,\|f_0^{(i)}\|_1\,\Bigl\|\Bigl(\sum_{k\neq0}
    \theta_k\varphi_k\Bigr)^{(\beta-i)}\Bigr\|_2+
    \Bigl\|\sum_{i=\lceil\beta\rceil}^{\infty} {\beta\choose i}
    \ (f_0\ast g_A)^{(i)}\Bigl(\sum_{k\neq0}
    \theta_k\varphi_k\Bigr)^{(\beta-i)}\Bigr\|_2\\\label{37}
&=&\sum_{i=0}^{\lfloor\beta\rfloor} {\beta\choose i}
    \ \frac{\|f_0^{(i)}\|_1}{2A-1}\,\Bigl\|\Bigl(\sum_{k\neq0}
    \theta_k\varphi_k\Bigr)^{(\beta-i)}\Bigr\|_2+
    \Bigl\|\sum_{i=\lceil\beta\rceil}^{\infty} {\beta\choose i}
    \ (f_0\ast g_A)^{(i)}\Bigl(\sum_{k\neq0}
    \theta_k\varphi_k\Bigr)^{(\beta-i)}\Bigr\|_2
\end{eqnarray}
where $\|f_0^{(i)}\|_1$ is of course equal to 1 for $i=0$ and
finite for $i=1,\ldots, \lfloor\beta\rfloor$. When
$\beta\in\mathbb N$, then ${\beta\choose i}=0$ for all
$i\ge\lceil\beta\rceil$, so there is no residual. In the next step
we employ:
\begin{eqnarray}
\label{32}
\Bigl\|\Bigl(\sum_{k\neq0}\theta_k\varphi_k\Bigr)^{(\gamma)}\Bigr\|_2^2
&=&\sum_{k\neq0} \theta_k^2\left(\frac{k\pi}{A}\right)^{2\gamma}
\end{eqnarray}
for all $\gamma$, proven in (\ref{34}). Furthermore for
$\Theta_A(L)$, $\sum|\theta_k|\le A^{-2\beta+1}$ and
$\sum\theta_k^2\left(\frac{k\pi}{A}\right)^{2\beta}\le 4A^2L$ had
been determined in (\ref{8}). Therefrom we can show in (\ref{35})
that
\begin{eqnarray}
\label{33}
\sum_{k\neq0}\theta_k^2\left(\frac{k\pi}{A}\right)^{2\beta-l}&\le&
(1+4L)A^{2-l}\qquad\mbox{for}\ \ 0<\l<2\beta
\end{eqnarray}
Hence continuing at inequality number (\ref{37}):
\begin{eqnarray}\nonumber
&&\Bigl\|\Bigl(f_0\ast g_A\cdot\sum_{k\neq0}
    \theta_k\varphi_k\Bigr)^{(\beta)}\Bigr\|_2\\\nonumber
    \nonumber
&=&\sum_{i=0}^{\lfloor\beta\rfloor} {\beta\choose i}
    \ \frac{\|f_0^{(i)}\|_1}{2A-1}\,\sqrt{\sum_{k\neq0}
    \theta_k^2\left(\frac{k\pi}{A}\right)^{2(\beta-i)}}+
    \Bigl\|\sum_{i=\lceil\beta\rceil}^{\infty} {\beta\choose i}
    \ (f_0\ast g_A)^{(i)}\Bigl(\sum_{k\neq0}
    \theta_k\varphi_k\Bigr)^{(\beta-i)}\Bigr\|_2\\\nonumber
&\le&\frac{1}{2A-1}\,\sqrt{4LA^2} +
    \sum_{i=1}^{\lfloor\beta\rfloor} {\beta\choose i}
    \ \frac{\|f_0^{(i)}\|_1}{2A-1}\,\sqrt{(1+4L)A^{2-i}}\\\label{38}
&&+\ \Bigl\|\sum_{i=\lceil\beta\rceil}^{\infty} {\beta\choose i}
    \ (f_0\ast g_A)^{(i)}\Bigl(\sum_{k\neq0}
    \theta_k\varphi_k\Bigr)^{(\beta-i)}\Bigr\|_2
\end{eqnarray}
For the residual we apply Lemma 3 to our functions. It states that
for functions with support in $[-A,A]$:
\begin{eqnarray*}
\Bigl\|\sum_{i=\lceil\beta\rceil}^{\infty} {\beta\choose i}
    \ f^{(i)}\cdot g^{(\beta-i)}\Bigr\|_2^2&=&
    o(A^2)\ \|\widehat{f
    ^{(\lceil\beta\rceil)}}\|_1^2\
    \|g\|_2^2
\end{eqnarray*}
Setting $f:=f_0\ast g_A$ and $g:=\sum\theta_k\varphi_k$, we
proceed at inequality number (\ref{38}):
\begin{eqnarray}\nonumber
\Bigl\|\Bigl(f_0\ast g_A\cdot\sum_{k\neq0}
    \theta_k\varphi_k\Bigr)^{(\beta)}\Bigr\|_2
    &\le&\frac{1}{2A-1}\,\sqrt{4LA^2} +
    \sum_{i=1}^{\lfloor\beta\rfloor} {\beta\choose i}
    \ \frac{\|f_0^{(i)}\|_1}{2A-1}\,\sqrt{(1+4L)A^{2-i}}\\\nonumber
&&+\ o(A)\ \Bigl\|\widehat{(f_0\ast g_A)
    ^{(\lceil\beta\rceil)}}\Bigr\|_1\
    \Bigl\|\sum_{k\neq0}\theta_k\varphi_k\Bigr\|_2\\\nonumber
&=&\frac{1}{2A-1}\,\sqrt{4LA^2} +
    \sum_{i=1}^{\lfloor\beta\rfloor} {\beta\choose i}
    \ \frac{\|f_0^{(i)}\|_1}{2A-1}\,\sqrt{(1+4L)A^{2-i}}\\\label{39}
&&+\ o(A)\ \Bigl\|\widehat{(f_0\ast g_A)
    ^{(\lceil\beta\rceil)}}\Bigr\|_1\
    \sqrt{\sum_{k\neq0}\theta_k^2}
\end{eqnarray}
After having derived the claim of Theorem 2, we will show in
(\ref{42}) that
\begin{eqnarray}
\label{40} \|\widehat{(f_0\ast
    g_A)^{(\lceil\beta\rceil)}}\|_1\ \le\ \frac{\|\widehat{f_0^{(\lfloor
    \beta\rfloor)}}\|_1}{A-1/2},&& \mbox{where}\quad
    \smash{\|\widehat{f_0^{(\lfloor \beta\rfloor)}}\|_1\ <\ \infty}
\end{eqnarray}
Furthermore $\sum\theta_k^2\le\sum|\theta_k|\le A^{-2\beta+1}$,
but $-2\beta+1<0$, such that (\ref{39}) can be continued as
\begin{eqnarray}\nonumber
\Bigl\|\Bigl(f_0\ast g_A\cdot\sum_{k\neq0}
    \theta_k\varphi_k\Bigr)^{(\beta)}\Bigr\|_2
    &<&\frac{1}{2A-1}\,\sqrt{4LA^2} +
    \sum_{i=1}^{\lfloor\beta\rfloor} {\beta\choose i}
    \ \frac{\|f_0^{(i)}\|_1}{2A-1}\, \sqrt{(1+4L)A^{2-i}}\\\nonumber
&&+\ o(A)\ \frac{\|\widehat{f_0
    ^{(\lfloor\beta\rfloor)}}\|_1}{A-1/2}\
    \sqrt{A^{-2\beta+1}}\\\nonumber
&=&\frac{1}{2A-1}\,\sqrt{4LA^2} +
    \sum_{i=1}^{\lfloor\beta\rfloor} {\beta\choose i}
    \ \frac{\|f_0^{(i)}\|_1}{2A-1}\, \sqrt{(1+4L)A^{2-i}}\\\nonumber
&&+\ o(A)\
    \frac{\frac{1}{2\pi}\|\widehat{f_0
    ^{(\lfloor\beta\rfloor)}}\|_1}{A-1/2}\
    o(1)\\\nonumber
&=&\sqrt{L}\Bigl(1+o(1)\Bigr)+O\left(A^{-1/2}\right)+o(1)
\end{eqnarray}
This result in connection with (\ref{21}) and (\ref{36}) completes
Theorem 2:
\begin{eqnarray*}
\qquad\qquad\|f_{\theta}^{(\beta)}\|_2&\le&\frac{1}{b(\theta)}\Bigl\|\left(f_0\ast
    g_A\right)^{(\beta)}\Bigr\|_2+
    \frac{1}{b(\theta)}\Bigl\|\Bigl(f_0\ast g_A\sum\limits_{k\neq0} \theta_k\varphi_k\Bigr)
    ^{(\beta)}\Bigr\|_2\\
&=&O\left(A^{-1}\right)\|f_0^{(\beta-1)}\|_2
    +\sqrt{L}\Bigl(1+o(1)\Bigr)\\
&=&\sqrt{L}+o(1)
\end{eqnarray*}
Still we are left to prove the intermediate assertions (\ref{32}),
(\ref{33}) and (\ref{40}).

As an exception to the ordinary case, sine and cosine enjoy an
easy to calculate fractional derivative:
$\sin^{(\gamma)}(ax)=a^{\gamma}\sin\ (ax+\gamma\pi/2)$ and the
like for cosine (Samko \cite{sam}, p. 174). Obviously, the
orthogonality between our functions $\varphi_k$ is preserved
through derivation.
\begin{eqnarray}\nonumber
&&\int\Bigl(\sum_{k\neq0}\theta_k\varphi_k^{(\gamma)}(x)\Bigr)^2dx\\\nonumber
&=&\int\sum_{k\neq0}\theta_k^2\
    \varphi_k^{(\gamma)}(x)^2dx\\\nonumber
&=&\int_{-A}^A
    \sum_{k>0}\theta_k^2\ \frac{1}{A}
    \left(\frac{k\pi}{A}\right)^{2\gamma}
    \cos^2\left(\frac{k\pi x}{A}+\frac{\gamma\pi}{2}\right)+
    \sum_{k<0}\theta_k^2\ \frac{1}{A}
    \left(\frac{k\pi}{A}\right)^{2\gamma}
    \sin^2\left(\frac{k\pi x}{A}+\frac{\gamma\pi}{2}\right)\ dx\\\nonumber
&=&\sum_{k>0}\theta_k^2\
    \left(\frac{k\pi}{A}\right)^{2\gamma}
    \frac{1}{A}\int_{-A}^A\cos^2\left(\frac{k\pi x}{A}+\frac{\gamma\pi}{2}\right)+
    \sum_{k<0}\theta_k^2
    \left(\frac{k\pi}{A}\right)^{2\gamma}
    \frac{1}{A}\int_{-A}^A\sin^2\left(\frac{k\pi x}{A}+\frac{\gamma\pi}{2}\right)\ dx\\\nonumber
&=&\sum_{k>0}\theta_k^2\
    \left(\frac{k\pi}{A}\right)^{2\gamma}
    \frac{1}{A}\int_{-A}^A\cos^2\frac{k\pi x}{A}\ dx+
    \sum_{k<0}\theta_k^2
    \left(\frac{k\pi}{A}\right)^{2\gamma}
    \frac{1}{A}\int_{-A}^A\sin^2\frac{k\pi x}{A}\ dx\\\nonumber
&=&\sum_{k>0}\theta_k^2\
    \left(\frac{k\pi}{A}\right)^{2\gamma}
    \frac{1}{A}\ A+
    \sum_{k<0}\theta_k^2
    \left(\frac{k\pi}{A}\right)^{2\gamma}
    \frac{1}{A}\ A\\
    \label{34}
&=&\sum_{k\neq0}\theta_k^2\
    \left(\frac{k\pi}{A}\right)^{2\gamma}
\end{eqnarray}
Referring to step (\ref{33}), $0<l<2\beta$:
\begin{eqnarray}\nonumber
\sum_{k\neq0}\theta_k^2\left(\frac{k\pi}{A}\right)^{2\beta-l}
&=&\sum_{0\neq|k|\le
    A^2/\pi}\theta_k^2\left(\frac{k\pi}{A}\right)^{2\beta-l}+
    \sum_{|k|>A^2/\pi}\theta_k^2\left(\frac{k\pi}{A}\right)^{2\beta}\left(\frac{k\pi}{A}\right)^{-l}\\\nonumber
&\le&A^{2\beta-l}\sum_{0\neq|k|\le A^2/\pi}\theta_k^2+
    A^{-l}\sum_{|k|>A^2/\pi}\theta_k^2\left(\frac{k\pi}{A}\right)^{2\beta}\\\nonumber
&\le&A^{2\beta-l}\sum_{k\neq0}|\theta_k|+
    A^{-l}\sum_{k\neq0}\theta_k^2\left(\frac{k\pi}{A}\right)^{2\beta}\\\nonumber
&\le&A^{2\beta-l}\cdot A^{-2\beta+1}+
    A^{-l}\cdot4LA^2\\\nonumber
&=&A^{1-l}+4LA^{2-l}\\
    \label{35}
&\le&(1+4L)A^{2-l}
\end{eqnarray}
Proof of (\ref{40}):
\begin{eqnarray}\nonumber
\|\widehat{(f_0\ast g_A)^{(\lceil\beta\rceil)}}\|_1
    &=&\int \Bigl|\omega^{\lceil\beta\rceil} \widehat f_0(\omega)
    \widehat g_A(\omega)\Bigr|d\omega\\\nonumber
&=&\int \Bigl|\omega^{\lceil\beta\rceil} \widehat f_0(\omega)\,
    \frac{2\sin(A-1/2)\omega}{(2A-1)\omega}\Bigr|d\omega\\\nonumber
&=&\frac{1}{A-1/2}\,\int \Bigl|\omega^{\lfloor\beta\rfloor}
    \widehat f_0(\omega)\,
    \sin(A-1/2)\omega\Bigr|d\omega\\\label{42}
&\le& \frac{\|\widehat{f_0^{(\lfloor
    \beta\rfloor)}}\|_1}{A-1/2}
\end{eqnarray}
$\|\widehat{f_0^{(\lfloor\beta\rfloor)}}\|_1$ exists, because we
chose $f_0\in\mathcal S$. This concludes the proof of Theorem 2. \
\ \qquad\qquad$\square$
\\ \\
{\sl\textbf{Lemma 3}\quad For functions $f$ and $g$, which are
both $L_2$-integrable, sufficiently regular and have support in
$[-A,A]$, it holds that}
\begin{eqnarray*}
\Bigl\|\sum_{i=\lceil\beta\rceil}^{\infty}{\beta\choose i}
    \ f^{(i)} \cdot g^{(\beta-i)}\Bigr\|_2^2&=&
    o(A^2)\
    \Bigl\|\widehat {f^{\lceil\beta\rceil}}\Bigr\|_1^2\cdot\|g\|_2^2
\end{eqnarray*}
\textbf{Proof}\quad This proof takes a detour via Fourier
coefficients. Begin with the following discussion: The power
function is an analytical function. We may thus for instance
expand $(\frac{\kappa\pi}{A})^{\beta}$ into an infinite Taylor
series at point $\frac{(\kappa-\lambda)\pi}{A}$.
\begin{eqnarray*}
\Bigl(\frac{\kappa\pi}{A}\Bigr)^{\beta}
    &=&\sum_{i=0}^{\infty}{\beta\choose i}\
    \Bigl(\frac{\lambda\pi}{A}\Bigr)^i\,
    \Bigl(\frac{(\kappa-\lambda)\pi}{A}\Bigr)^{\beta-i}
\end{eqnarray*}
We cut the Tailor expansion of
$\Bigl(\frac{\kappa\pi}{A}\Bigr)^{\beta}$ after
$\lfloor\beta\rfloor$ and bound the residual.
\begin{eqnarray*}
&&\Bigl|\sum_{i=\lceil\beta\rceil}^{\infty}{\beta\choose i}\
    \Bigl(\frac{\lambda\pi}{A}\Bigr)^i\,
    \Bigl(\frac{(\kappa-\lambda)\pi}{A}\Bigr)^{\beta-i}\Bigr|   \\
&=&\Bigl|\sum_{i=0}^{\infty}{\beta\choose \lceil\beta\rceil+i}\
    \Bigl(\frac{\lambda\pi}{A}\Bigr)
    ^{\lceil\beta\rceil+i}\Bigl(\frac{(\kappa-\lambda)\pi}{A}\Bigr)
    ^{\beta-\lceil\beta\rceil-i}\Bigr|\\
&=&\Bigl|\sum_{i=0}^{\infty}\frac{
    \beta\cdots(\beta-\lceil\beta\rceil-i+1)}
    {(\lceil\beta\rceil+i)!}
    \Bigl(\frac{\lambda\pi}{A}\Bigr)
    ^{\lceil\beta\rceil+i}\Bigl(\frac{(\kappa-\lambda)\pi}{A}\Bigr)
    ^{\beta-\lceil\beta\rceil-i}\Bigr|\\
&=&\Bigl|\frac{\beta\cdots(\beta-\lceil\beta\rceil+1)}
    {\lceil\beta\rceil!}\ \Bigl(\frac{\lambda\pi}{A}\Bigr)^{\lceil\beta\rceil}
    \sum_{i=0}^{\infty}\frac{(\beta-\lceil\beta\rceil)
    \cdots(\beta-\lceil\beta\rceil-i+1)}
    {(\lceil\beta\rceil+i)\cdots(\lceil\beta\rceil+1)}\
    \Bigl(\frac{\lambda\pi}{A}\Bigr)^{i}\Bigl(\frac{(\kappa-\lambda)\pi}{A}\Bigr)
    ^{\beta-\lceil\beta\rceil-i}\Bigr|\\
&\le&{\beta\choose\lceil\beta\rceil}\
    \Bigl|\frac{\lambda\pi}{A}\Bigr|^{\lceil\beta\rceil}
    \sum_{i=0}^{\infty}\ \Bigl|\frac{(\beta-\lceil\beta\rceil)
    \cdots(\beta-\lceil\beta\rceil-i+1)}
    {i!}\
    \Bigl(\frac{\lambda\pi}{A}\Bigr)^{i}\Bigl(\frac{(\kappa-\lambda)\pi}{A}\Bigr)
    ^{\beta-\lceil\beta\rceil-i}\Bigr|
\end{eqnarray*}
The product $(\beta-\lceil\beta\rceil)
\cdots(\beta-\lceil\beta\rceil-i+1)$ consists of $i$ factors,
which are all negative. We can write $|(\beta-\lceil\beta\rceil)
\cdots(\beta-\lceil\beta\rceil-i+1)| =
(-1)^i(\beta-\lceil\beta\rceil)
\cdots(\beta-\lceil\beta\rceil-i+1)$, such that
\begin{eqnarray*}
&&\Bigl|\sum_{i=\lceil\beta\rceil}^{\infty}{\beta\choose i}\
    \Bigl(\frac{\lambda\pi}{A}\Bigr)^i\,
    \Bigl(\frac{(\kappa-\lambda)\pi}{A}\Bigr)^{\beta-i}\Bigr|   \\
&\le&{\beta\choose\lceil\beta\rceil}\
    \Bigl|\frac{\lambda\pi}{A}\Bigr|^{\lceil\beta\rceil}
    \sum_{i=0}^{\infty}\ (-1)^i\frac{(\beta-\lceil\beta\rceil)
    \cdots(\beta-\lceil\beta\rceil-i+1)}
    {i!}\
    \Bigl|\frac{\lambda\pi}{A}\Bigr|^{i}\Bigl|\frac{(\kappa-\lambda)\pi}{A}\Bigr|
    ^{\beta-\lceil\beta\rceil-i}\\
&=&{\beta\choose\lceil\beta\rceil}\
    \Bigl|\frac{\lambda\pi}{A}\Bigr|^{\lceil\beta\rceil}
    \sum_{i=0}^{\infty}{\beta-\lceil\beta\rceil\choose i}\
    \Bigl(\frac{-|\lambda|\pi}{A}\Bigr)^{i}\Bigl(\frac{|\kappa-\lambda|\pi}{A}\Bigr)
    ^{\beta-\lceil\beta\rceil-i}\\
&=&{\beta\choose\lceil\beta\rceil}\
    \Bigl|\frac{\lambda\pi}{A}\Bigr|^{\lceil\beta\rceil}
    \Bigl(\frac{(|\kappa-\lambda|-|\lambda|)\pi}{A}\Bigr)
    ^{\beta-\lceil\beta\rceil}
\end{eqnarray*}
Since we know that $-1<\beta-\lceil\beta\rceil<0$, we can
approximte
$(\frac{(|\kappa-\lambda|-|\lambda|)\pi}{A})^{\beta-\lceil\beta\rceil}
=O(A^{\beta-\lceil\beta\rceil})=o(A)$. Now we expand the tail of
our Leibnitz formula into a Fourier series and plug in the bound
of the Taylor series:
\begin{eqnarray}\nonumber
&&\Bigl\|\sum_{i=\lceil\beta\rceil}^{\infty}{\beta\choose i}
    \ f^{(i)} \cdot g^{(\beta-i)}\Bigr\|_2^2\\
    \nonumber
&=&\frac{1}{2A}
    \sum_{\kappa\in\mathbb Z}\left(\sum_{i=\lceil\beta\rceil}^
    {\infty}{\beta\choose i}\ \widehat{f^{(i)}} \ast
    \widehat{g^{(\beta-i)}}\Bigl(\frac{\kappa\pi}{A}\Bigr)\right)^2\\
    \nonumber
&=&\frac{1}{2A}
    \sum_{\kappa\in\mathbb Z}\left(\sum_{i=\lceil\beta\rceil}^
    {\infty}{\beta\choose i}\ \frac{1}{2A}\sum_{\lambda\in\mathbb Z}
    \widehat{f^{(i)}}\Bigl(\frac{\lambda\pi}{A}\Bigr)
    \widehat{g^{(\beta-i)}}\Bigl(\frac{(\kappa-\lambda)\pi}{A}\Bigr)\right)^2\\
    \nonumber
&=&\frac{1}{2A}
    \sum_{\kappa\in\mathbb Z}\left(\sum_{i=\lceil\beta\rceil}^
    {\infty}{\beta\choose i}\ \frac{1}{2A}\sum_{\lambda\in\mathbb Z}
    \Bigl(\frac{\lambda\pi}{A}\Bigr)^i\widehat{f}\Bigl(\frac{\lambda\pi}{A}\Bigr)
    \Bigl(\frac{(\kappa-\lambda)\pi}{A}\Bigr)^{\beta-i}
    \widehat{g}\Bigl(\frac{(\kappa-\lambda)\pi}{A}\Bigr)\right)^2\\
    \nonumber
&=&\frac{1}{2A}
    \sum_{\kappa\in\mathbb Z}\left(\frac{1}{2A}\sum_{\lambda\in\mathbb Z}
    \sum_{i=\lceil\beta\rceil}^{\infty}{\beta\choose i}\
    \Bigl(\frac{\lambda\pi}{A}\Bigr)^i \Bigl(\frac{(\kappa-\lambda)\pi}{A}\Bigr)^{\beta-i}
    \widehat{f}\Bigl(\frac{\lambda\pi}{A}\Bigr)\
    \widehat{g}\Bigl(\frac{(\kappa-\lambda)\pi}{A}\Bigr)\right)^2\\
    \nonumber
&\le&\frac{1}{2A}
    \sum_{\kappa\in\mathbb Z}\left(\frac{1}{2A}\sum_{\lambda\in\mathbb Z}
    o(A){\beta\choose \lceil\beta\rceil}\
    \Bigl|\frac{\lambda\pi}{A}\Bigr|^{\lceil\beta\rceil}
    \Bigl|\widehat{f}\Bigl(\frac{\lambda\pi}{A}\Bigr)\Bigr|\cdot
    \Bigl|\widehat{g}\Bigl(\frac{(\kappa-\lambda)\pi}{A}\Bigr)\Bigr|\right)^2\\
    \nonumber
&=&o(A^2)\ \frac{1}{2A}
    \sum_{\kappa\in\mathbb Z}\frac{1}{2A}\sum_{\lambda\in\mathbb Z}
    {\beta\choose \lceil\beta\rceil}\
    \Bigl|\frac{\lambda\pi}{A}\Bigr|^{\lceil\beta\rceil}
    \Bigl|\widehat{f}\Bigl(\frac{\lambda\pi}{A}\Bigr)\Bigr|\cdot
    \Bigl|\widehat{g}\Bigl(\frac{(\kappa-\lambda)\pi}{A}\Bigr)\Bigr|\\
    \nonumber
&&\qquad\quad\times\ \frac{1}{2A}\sum_{\mu\in\mathbb Z}
    {\beta\choose \lceil\beta\rceil}\
    \Bigl|\frac{\mu\pi}{A}\Bigr|^{\lceil\beta\rceil}
    \Bigl|\widehat{f}\Bigl(\frac{\mu\pi}{A}\Bigr)\Bigr|\cdot
    \Bigl|\widehat{g}\Bigl(\frac{(\kappa-\mu)\pi}{A}\Bigr)\Bigr|\\
    \nonumber
&=&o(A^2)\ {\beta\choose \lceil\beta\rceil}^2\
    \frac{1}{2A}\sum_{\lambda\in\mathbb Z}
    \Bigl|\frac{\lambda\pi}{A}\Bigr|^{\lceil\beta\rceil}
    \Bigl|\widehat{f}\Bigl(\frac{\lambda\pi}{A}\Bigr)\Bigr|\ \
    \frac{1}{2A}\sum_{\mu\in\mathbb Z}\Bigl|\frac{\mu\pi}{A}\Bigr|^{\lceil\beta\rceil}
    \Bigl|\widehat{f}\Bigl(\frac{\mu\pi}{A}\Bigr)\Bigr|\\
    \nonumber
&&\qquad\quad\times\ \frac{1}{2A}\sum_{\kappa\in\mathbb Z}
    \Bigl|\widehat{g}\Bigl(\frac{(\kappa-\lambda)\pi}{A}\Bigr)\Bigr|\cdot
    \Bigl|\widehat{g}\Bigl(\frac{(\kappa-\mu)\pi}{A}\Bigr)\Bigr|\\
    \nonumber
&\le&o(A^2)\ {\beta\choose \lceil\beta\rceil}^2\
    \frac{1}{2A}\sum_{\lambda\in\mathbb Z}
    \Bigl|\Bigl(\frac{\lambda\pi}{A}\Bigr)^{\lceil\beta\rceil}
    \widehat{f}\Bigl(\frac{\lambda\pi}{A}\Bigr)\Bigr|\ \
    \frac{1}{2A}\sum_{\mu\in\mathbb Z}\Bigl|\Bigl(\frac{\mu\pi}{A}\Bigr)^{\lceil\beta\rceil}
    \widehat{f}\Bigl(\frac{\mu\pi}{A}\Bigr)\Bigr|\\
    \nonumber
&&\qquad\quad\times\ \sqrt{\frac{1}{2A}\sum_{\kappa\in\mathbb Z}
    \Bigl|\widehat{g}\Bigl(\frac{(\kappa-\lambda)\pi}{A}\Bigr)\Bigr|^2}
    \sqrt{\frac{1}{2A}\sum_{\kappa\in\mathbb Z}
    \Bigl|\widehat{g}\Bigl(\frac{(\kappa-\mu)\pi}{A}\Bigr)\Bigr|^2}\\
    \nonumber
&=&o(A^2)\ {\beta\choose \lceil\beta\rceil}^2\
    \left(\frac{1}{2A}\sum_{\lambda\in\mathbb Z}
    \Bigl|\Bigl(\frac{\lambda\pi}{A}\Bigr)^{\lceil\beta\rceil}
    \widehat{f}\Bigl(\frac{\lambda\pi}{A}\Bigr)\Bigr|\right)^2\cdot\|g\|_2^2
\end{eqnarray}
For growing $A$, the Fourier expansion approaches the Fourier
transform, and hence
\begin{eqnarray}
&&\Bigl\|\sum_{i=\lceil\beta\rceil}^{\infty}{\beta\choose i}
    \ f^{(i)} \cdot g^{(\beta-i)}\Bigr\|_2^2\\
    \nonumber
&=&o(A^2)\ {\beta\choose \lceil\beta\rceil}^2\
    \left(\frac{1}{2\pi}\int
    \Bigl|\omega^{\lceil\beta\rceil}
    \widehat{f}(\omega)\Bigr|d\omega\Bigl(1+o(1)\Bigr)\right)^2\cdot\|g\|_2^2\\
    \nonumber
&=&o(A^2)\ \Bigl\|\widehat
    {f^{\lceil\beta\rceil}}\Bigr\|_1^2\cdot\|g\|_2^2
\end{eqnarray}
which is the statement of Lemma 3.
\qquad\qquad\qquad\qquad\qquad\qquad\qquad\qquad
\qquad\qquad\qquad\quad$\square$
\\ \\

\textbf{Proof} of \textbf{Lemma 2}\quad We start with
\begin{eqnarray}\nonumber
&&I_{f_{\theta}}(\theta_{\kappa})\\
    \nonumber
&=&\int_{-A}^A\frac{\left(\partial
    f_{\theta}(x)/\partial
    \theta_{\kappa}\right)^2}{f_{\theta}(x)}\ dx\\\nonumber
&=&\frac{1}{b^2(\theta)}\int_{-A}^A\frac{1}{f_{\theta}(x)}\left[-f_{\theta}(x)\int_{-A}^A
    f_0\ast g_A(y)\varphi_{\kappa}(y)\ dy
    +f_0\ast g_A(x)\varphi_{\kappa}(x)\right]^2dx\\\nonumber
&=&\frac{1}{b^2(\theta)}\int_{-A}^A\left[
    f_{\theta}(x)\left(\int_{-A}^A f_0\ast
    g_A(y)\varphi_{\kappa}(y)\ dy\right)^2-2 \int_{-A}^A f_0\ast
    g_A(y)\varphi_{\kappa}(y)\ dy\ f_0\ast
    g_A(x)\varphi_{\kappa}(x)\right.\\\nonumber
&&\quad\left.+\ \frac{1}{f_{\theta}(x)}\Bigl(f_0\ast
    g_A(x)\varphi_{\kappa}(x)\Bigr)^2\right] dx\\
    \nonumber
&=&-\frac{1}{b^2(\theta)}\left[\int_{-A}^A f_0\ast
    g_A(y)\varphi_{\kappa}(y)\ dy\right]^2+
    \frac{1}{b^2(\theta)}\int_{-A}^A \frac{\Bigl(f_0\ast
    g_A(x)\Bigr)^2\varphi^2_{\kappa}(x)}
    {\frac{1}{b(\theta)}f_0\ast g_A(x)\Bigl(1+
    \sum\limits_{\lambda\neq0}
    \theta_{\lambda}\varphi_{\lambda}(x)\Bigr)} dx\\
    \nonumber
&=&-\frac{1}{b^2(\theta)}\left[\int_{-A}^A f_0\ast
    g_A(y)\varphi_{\kappa}(y)\ dy\right]^2+
    \frac{1}{b(\theta)}\int_{-A}^A \frac{f_0\ast g_A(x)\varphi_{\kappa}^2(x)}{1+
    \sum\limits_{\lambda\neq0}
    \theta_{\lambda}\varphi_{\lambda}(x)}\ dx\\
    \nonumber
&=&-\frac{1}{b^2(\theta)}\left[\int_{-A}^A\frac{\varphi_{\kappa}(x)}{2A-1}\
    dx-
    2\int_{A-1}^A\left(\frac{1}{2A-1}-f_0\ast g_A(x)\right)\varphi_{\kappa}(x)\ dx\right]^2\\
    \nonumber
&&+\ \frac{1}{b(\theta)}\left[\int_{-A}^A
    \frac{\frac{1}{2A-1}\ \varphi_{\kappa}^2(x)}
    {1+\sum\limits_{\lambda\neq0}
    \theta_{\lambda}\varphi_{\lambda}(x)} \ dx-
    2\int_{A-1}^A\frac{\left(\frac{1}{2A-1}-f_0\ast
    g_A(x)\right)\varphi_{\kappa}^2(x)}
    {1+\sum\limits_{\lambda\neq0}
    \theta_{\lambda}\varphi_{\lambda}(x)}\ dx\right]
\end{eqnarray}
The leading term is $\int_{-A}^A \frac{\frac{1}{2A-1}\
\varphi_{\kappa}^2(x)}{1+\sum\theta_{\lambda}\varphi_{\lambda}(x)}
\ dx$. Due to $|\sum\theta_{\lambda}\varphi_{\lambda} (x)|\le
A^{-1/2}\sum|\theta_{\lambda}|\le A^{-1/2}\cdot\nolinebreak
A^{-2\beta+1}$ $< A^{-1/2}$, we know it lies in the interval
$((2A-1)^{-1}(1+A^{-1/2})^{-1},(2A-1)^{-1}(1-A^{-1/2})^{-1})$.
Moreover from (\ref{21}) we have $1-A^{-3/2}\le
b(\theta)\le1+A^{-3/2}$. For $A\longrightarrow\infty$ we obtain:
\begin{eqnarray}
    \nonumber
I_{f_{\theta}}(\theta_{\kappa}) &=&\frac{1}{b^2(\theta)}\left[\
    0+O\left(\frac{1}{(A-1/2)\sqrt{A}}\right)\right]^2+
    \frac{1}{b(\theta)}\left[\frac{1+o(1)}{2A-1}+
    O\left(\frac{1+o(1)}{(A-1/2)A}\right)\right]\\
    \nonumber
&=&\Bigl(1+o(1)\Bigr)O\left(A^{-3}\right)+
    \Bigl(1+o(1)\Bigr)\left[\frac{1+o(1)}{2A-1}+
    O\left(A^{-2}\right)\right]\\
    \nonumber
&=&\frac{1+o(1)}{2A}
\end{eqnarray}
$\mbox{Re}\widehat f_{\theta}\left(\frac{{\kappa}\pi}{A}\right)$
can be expressed as $A^{1/2}\int
f_{\theta}(x)\varphi_{|\kappa|}(x)dx$, yielding
\begin{eqnarray}\nonumber
&&\frac{\partial\,\mbox{Re}\widehat
    f_{\theta}\left(\frac{{\kappa}\pi}{A}\right)}{\partial
    \theta_{|\kappa|}}\\
    \nonumber
&=&\frac{\partial}{\partial
    \theta_{|\kappa|}}\sqrt{A}\int_{-A}^A
    f_{\theta}(x)\varphi_{|\kappa|}(x)dx\\\nonumber
&=&\frac{\sqrt{A}}{b(\theta)}\left[-\int_{-A}^A
    f_{\theta}(x)\!
    \int_{-A}^A\! f_0\ast g_A(y)\varphi_{|\kappa|}(y)dy
    \ \varphi_{|\kappa|}(x)dx + \int_{-A}^A\! f_0\ast g_A(x)\varphi_
    {|\kappa|}^2(x)dx\right]\\\nonumber
&=&\frac{\sqrt{A}}{b(\theta)}\left[-\frac{\mbox{Re}\widehat
    f_{\theta}\left(\frac{{\kappa}\pi}{A}\right)}{\sqrt{A}}
    \int_{-A}^A\! f_0\ast g_A(y)\varphi_{|\kappa|}(y)dy
    + \int_{-A}^A\! f_0\ast g_A(x)\varphi_{|\kappa|}^2(x)dx\right]
    \\\nonumber
&=&\frac{\mbox{Re}\widehat
    f_{\theta}\left(\frac{{\kappa}\pi}{A}\right)}{b(\theta)}
    \left[-\int_{-A}^A\!\frac{1}{2A-1}\ \varphi_{|\kappa|}(y)dy+2
    \int_{A-1}^A\!\left(\frac{1}{2A-1}-f_0\ast g_A(y)\right)\varphi_{|\kappa|}(y)dy
    \right]\\\nonumber
&&+\ \frac{\sqrt{A}}{b(\theta)}\left[ \int_{-A}^A\!\frac{1}{2A-1}\
    \varphi_{|\kappa|}^2(x)dx-2
    \int_{A-1}^A\!\left(\frac{1}{2A-1}-f_0\ast g_A(x)\right)\varphi_{|\kappa|}^2(x)dx
    \right]\\\nonumber
&=&\frac{\mbox{Re}\widehat
    f_{\theta}\left(\frac{{\kappa}\pi}{A}\right)}{1+o(1)}
    \left[\ 0+O\left(\frac{1}{(A-1/2)\sqrt{A}}\right)\right]+\
    \frac{\sqrt{A}}{1+o(1)}\left[ \frac{1}{2A-1}+O
    \left(\frac{1}{(A-1/2)A}\right)\right]\\
    \nonumber
&=&\frac{1+o(1)}{2\sqrt{A}}
\end{eqnarray}
A similar result is obtained for $\mbox{Im}\widehat
f_{\theta}\left(\frac{{\kappa}\pi}{A}\right)$, whereby $\int
f_0\ast g_A(y)\varphi_{-|\kappa|}(y)dy=0$ simplifies the task,
because the sine function is
anti-symmetric.\qquad\qquad\qquad\qquad
\qquad\qquad\qquad\qquad\qquad $\square$
\\\\
\textbf{Acknowledgement:} I would like to thank Prof. M. Neumann
for initiating and supporting the present work.

\end{document}